\documentclass[twoside,12pt]{article}
\usepackage{times}
\usepackage{easyReview}

\usepackage{setspace}

\usepackage{url}
\usepackage{emptypage}
\usepackage{nameref}

\usepackage{amsmath,amstext,amsgen,amsbsy,amsopn}
\usepackage{amssymb,amsfonts}
\usepackage{amsthm}
\usepackage{latexsym,amsxtra,euscript,amscd}

\usepackage[numbers]{natbib} 
\usepackage{nicematrix}
\usepackage{arydshln}

\usepackage{float}
\usepackage{xcolor}

\usepackage{titlesec}
\titleformat*{\section}{\normalfont\Large\bfseries\color{blue}}
\titleformat*{\subsection}{\normalfont\large\bfseries\color{blue}}
\titleformat*{\subsubsection}{\normalfont\normalsize\bfseries\color{blue}}
\usepackage[font={color=blue,bf}]{caption}
\usepackage[colorlinks=true,urlcolor=blue,linkcolor=black]{hyperref}

\newtheoremstyle{mystyle}{3pt}{3pt}
{\itshape\color{black}}{}
{\bfseries\color{blue}}{.}{.5em}{}
\theoremstyle{mystyle}

\usepackage{graphics,graphicx}
\graphicspath{{fig/},{../fig/}{../Fig/},{Fig/},{Fig_recu/}}
\usepackage{epsf}

\usepackage{enumitem}

\newtheorem{theorem}{Theorem}

\newtheorem{remark}{Remark}

\newcommand{\Det}[1]{|\!| #1 \text{\reflectbox{$|\!|$}}}
\DeclareMathOperator{\tr}{tr}

\def\A{\mathbf{A}}

\def\C{\mathbb{C}}

\usepackage{booktabs}

\usepackage{nicematrix}

\usepackage{tikz}
\usetikzlibrary{arrows,chains,matrix,positioning,scopes}
\usetikzlibrary{positioning,shadows,backgrounds}  
\tikzset{>=stealth',every on chain/.append style={join},
	every join/.style={->}}
\tikzstyle{labeled}=[execute at begin node=$\scriptstyle,
execute at end node=$]

\newcommand*{\ADRnl}{ORCID: 0000-0001-9062-9280, University of Sopron, Institute of Basic Sciences, Departement of Mathematics,  Bajcsy-Zs. u. 4, Sopron, 9400, Hungary. \texttt{nemeth.laszlo@uni-sopron.hu}}
\newcommand*{\ADRszl}{ORICD: 0000-0002-4582-6100, J.~Selye University, Department of Mathematics, Kom\'arno, Slovakia  
	and 
	University of Sopron, Institute of Basic Sciences, Departement of Mathematics,  Bajcsy-Zs. u. 4, Sopron, 9400, Hungary. \texttt{szalay.laszlo@uni-sopron.hu}}

\newcommand*{\TIT}{\color{blue}Explicit solution of system of two higher-order recurrences}
\title{\bf \TIT}

\author{L\'aszl\'o N\'emeth\footnote{\ADRnl} \  and L\'aszl\'o Szalay\footnote{\ADRszl} 
}

\begin{document}
	
	\maketitle \thispagestyle{empty}
	
	\begin{abstract}
		We give a method to determine an explicit solution to a system of two inhomogeneous linear recursive sequences of higher order. Our approach can be used efficiently in solving certain combinatorial problems. We finish the paper by considering a tiling problem with black and white dominoes, and we use the method as a demonstration to find the solution. 
		\\[1mm]
		{\em Key Words: linear recurrence, vector recurrence, tiling, $k$-generalized Fibonacci number.}\\
		{\em MSC code:  11B37, 05B45, 52C20.}\\ 
	\end{abstract}
	

	\section{Introduction}
	
	Recurrence sequences appear in several combinatorial problems. Time by time, there are more recurrences in the solution that possess the same recurrence relation. In this case, a suitable vector recurrence can simplify the description and the solution of the problem. Faye, N\'emeth, and Szalay~\cite{FNSz} investigated the question of 
	linear vector recurrences in general. Now we introduce some basic notations in order to formulate the problem.
	
	Let $s\ge1$ and $k\ge2$ denote two positive integers. Assume that there are given the matrices $\A_t=[a_{i,j}^{(t)}]\in\C^{k\times k}$ for $t=1,2,\dots,s$. We define the vector recurrence
	\begin{equation}\label{vecrec}
		\mathbf{v}_n=\mathbf{A}_1\mathbf{v}_{n-1}+\mathbf{A}_2\mathbf{v}_{n-2}+\cdots +\mathbf{A}_s\mathbf{v}_{n-s},\qquad n\ge s
	\end{equation}
	with initial column vectors $\mathbf{v}_{0},\mathbf{v}_{1},\dots,\mathbf{v}_{s-1}\in\C^k$.
	
	In \cite{FNSz} the authors developed a procedure to separate the component sequences $(v_i^{(t)})_{t\ge0}$ for $i=1,2,\dots,k$, and to give their own common recursive relation. The crucial point was to find the characteristic polynomial of a matrix derived from the initial conditions. If one finds the zeros of the characteristic polynomial, then the zeros, together with the initial values, can provide an explicit form for each component sequence. The difficulty that may arise in practice is finding the zeros precisely. The situation changes when we have restrictions on (\ref{vecrec}). In particular, we are able to gain a more precise final result if the size of matrices appearing in (\ref{vecrec}) is $2\times2$. This paper studies the vector recurrences with coefficient matrices $\A_t\in\C^{2\times 2}$ ($t=1,2,\dots,s$).
	
	Clearly, the simplest version of (\ref{vecrec}) is
	\begin{equation}\label{vr0}
		\mathbf{v}_{n}=\A_1\mathbf{v}_{n-1},
	\end{equation}
	when we have only one matrix, i.e., the order of vector recurrence is 1. Hence, the combination of component recurrences goes back only to one term. In this case, the characteristic polynomial of the matrix $\A_1$ implies a recurrence rule of order $k$ valid for any of the se\-parated component sequences of the vectors $\mathbf{v}_{n}$ (see, for example, \cite[Lemma 2.1]{NSz_Power}, or \cite{NL_grow}). (Note that the basic field in \cite[Lemma 2.1]{NSz_Power} is $\mathbb{R}$, but the statement can be extended without changes for the field $\mathbb{C}$ as well.)
	For non-homogeneous version of (\ref{vr0}), see  \cite{FNSz} (with two sequences) and \cite{NSz_Power} (generally).

	If we have the matrices $\A_i\in\C^{2\times 2}$ in (\ref{vecrec}), then (\ref{vecrec}) can be also given in the form 
	\begin{equation}\label{eqsys:def}
		\begin{array}{lclclclcl} 
			a_{n} & = & \alpha_{1,1}^{(1)} a_{n-1}+ \alpha_{1,2}^{(1)} b_{n-1} &+& \alpha_{1,1}^{(2)} a_{n-2}+ \alpha_{1,2}^{(2)} b_{n-2} +& \cdots &+& \alpha_{1,1}^{(s)} a_{n-s}+ \alpha_{1,2}^{(s)} b_{n-s}, \\ 
			b_{n} & = & \alpha_{2,1}^{(1)} a_{n-1}+ \alpha_{2,2}^{(1)} b_{n-1} &+& \alpha_{2,1}^{(2)} a_{n-2}+ \alpha_{2,2}^{(2)} b_{n-2} + &\cdots &+& \alpha_{2,1}^{(s)} a_{n-s}+ \alpha_{2,2}^{(s)} b_{n-s},
		\end{array} 
	\end{equation}
	where $a_{i}$, $b_{i}$ $(0\leq i \leq s-1)$ are the initial values.
	Clearly, now 
	\[\mathbf{v}_i=\begin{pmatrix} a_{i}\\ b_{i} \end{pmatrix}\quad  {\rm and}\quad
	\mathbf{A}_t=\begin{pmatrix}  \alpha_{1,1}^{(t)} & \alpha_{1,2}^{(t)}\\  \alpha_{2,1}^{(t)} & \alpha_{2,2}^{(t)} \end{pmatrix}
	\]
	for $0\leq i$ and $1\leq t\leq s$. 
	
	We will find explicitly the solution of \eqref{eqsys:def} for the recurrence sequences $(a_n)_{n=0}^\infty$ and $(b_n)_{n=0}^\infty$ in the  form  
	\begin{equation} \label{eq:z-K}
		z_n= c^{(s)}_{1}z_{n-1}+c^{(s)}_{2}z_{n-2}+\cdots+ c^{(s)}_{k}z_{n-k} = \sum_{i=1}^{k} c^{(s)}_{i}z_{n-i},
		 \quad n\geq k
	\end{equation}
	with suitable coefficients $c_{i}^{(s)}$ $(1\leq i\leq k)$.
	Here the usage of upper index $(s)$ of $c^{(s)}_{i}$ ($s$ is the order of (\ref{vecrec})) does not seem justified, but it becomes useful later when we find connection between different systems (see Theorem \ref{th:main}). Obviously, \eqref{eq:z-K} is equivalent to 
	\begin{equation} \label{eq:z-K_vec}
		z_n= (\mathbf{c}^{(s)})^T \mathbf{z},
	\end{equation}
	where 
	$(\mathbf{c}^{(s)})^T=\begin{pmatrix} c_1^{(s)}, c_2^{(s)}, \ldots, c_k^{(s)}\end{pmatrix}$ 
	and
	$\mathbf{z}^T=\begin{pmatrix} z_{n-1}, z_{n-2}, \ldots, z_{n-k}\end{pmatrix}$.
	
	We even define additional matrices as follows. Let  $\mathbf{A}_{i,j}$ ($1\leq i,j \leq s$) be a linear combination of $\mathbf{A}_i$ and $\mathbf{A}_j$ given by
	
	\[ \mathbf{A}_{i,j}=\mathbf{A}_{i} \begin{pmatrix}1&0\\0&0\end{pmatrix}  +  \mathbf{A}_{j} \begin{pmatrix}0&0\\0&1\end{pmatrix} = 
	\begin{pmatrix}  \alpha_{1,1}^{(i)} & \alpha_{1,2}^{(j)}\\  \alpha_{2,1}^{(i)} & \alpha_{2,2}^{(j)} \end{pmatrix}.\]  
	Furthermore, we define the function $\Det{.}$ as the sum 
	\[ \Det{\mathbf{A}^{(i,j)}}=|{\mathbf{A}_{i,j}}|+|{\mathbf{A}_{j,i}}|\]
	of determinants of two associated matrices.  
	Trivially, 
	$\Det{\mathbf{A}^{(i,j)}}=\Det{\mathbf{A}^{(j,i)}}$ and $\Det{\mathbf{A}^{(i,i)}}=2|\mathbf{A}_{i}|$. Finally, as usual, let $\tr(\mathbf{A}_i)$ denote the trace of matrix $\mathbf{A}_i$.
	
	Now we are ready to present the  main result of this paper.
	\begin{theorem} [Main theorem]\label{th:main}
		Sequences $(a_n)_{n=0}^\infty$ and $(b_n)_{n=0}^\infty$ defined in \eqref{eqsys:def} satisfy the homogeneous linear recurrence relation 
		\begin{equation}\label{eq:z}
			z_n= \sum_{i=1}^{2s} c^{(s)}_{i}z_{n-i}, \quad n\geq 2s,
    	\end{equation} 
		of order $k=2s$, or equivalently  
		\begin{equation*}
			z_n = (\mathbf{c}^{(s)})^T \mathbf{z},
		\end{equation*}
		where 
		$(\mathbf{c}^{(s)})^T=\begin{pmatrix}  c^{(s)}_{1},  c^{(s)}_{2},  \ldots,  c^{(s)}_{2s} \end{pmatrix}$
		and
		$\mathbf{z}^T=\begin{pmatrix}  z_{n-1}, z_{n-2},  \ldots,  z_{n-2s} \end{pmatrix}$.
		The coefficient vector $\mathbf{c}^{(s)}$ ($s\ge2$) can be  recursively given by 
		\begin{equation}\label{eqvect_recu}
			\mathbf{c}^{(s)}=
			\begin{pmatrix}\\ \\ \\ \\ \mathbf{c}^{(s-1)}\\ \\ \\0\\0\end{pmatrix}+
			\begin{pmatrix} 0\\ \vdots\\ 0\\
				\tr(\mathbf{A}_{s})\\
				-\Det{\mathbf{A}^{(1,s)}} \\
				\vdots\\
				-\Det{\mathbf{A}^{(s-1,s)}} \\
				0
			\end{pmatrix}+
			\begin{pmatrix} 0\\ \vdots\\ 0\\
				0\\  0\\  \vdots\\ 0 \\ -|\mathbf{A}_{s}|
			\end{pmatrix}
		\end{equation}  
		with the initial vector  
		\begin{equation*} \mathbf{c}^{(1)}=\begin{pmatrix} \tr(\mathbf{A}_1)\\ -|\mathbf{A}_1| \end{pmatrix}.
		\end{equation*}    
		
	\end{theorem} 
	
	\section{Proof of the main theorem}    
	
	\begin{proof}
		First we rewrite \eqref{eqsys:def} into a system of $2s$ homogeneous linear recurrences of order $1$ applying the substitutions $a_{i}^{(j)}:=a_{i-j}$ and $b_{i}^{(j)}:=b_{i-j}$ for $0\le j\le s-1$ and $1\leq i$ in the following way:
		\begin{equation*}\label{eqsys:ext}
			\begin{array}{lclclclcl} 
				a_{n}^{(0)} & = & \alpha_{1,1}^{(1)} a_{n-1}^{(0)}+ \alpha_{1,2}^{(1)} b_{n-1}^{(0)} &+& \alpha_{1,1}^{(2)} a_{n-1}^{(1)}+ \alpha_{1,2}^{(2)} b_{n-1}^{(1)} + &\cdots &+& \alpha_{1,1}^{(s)} a_{n-1}^{(s-1)}+ \alpha_{1,2}^{(s)} b_{n-1}^{(s-1)}, \\ 
				b_{n}^{(0)}  & = & \alpha_{2,1}^{(1)} a_{n-1}^{(0)}+ \alpha_{2,2}^{(1)} b_{n-1}^{(0)} &+& \alpha_{2,1}^{(2)} a_{n-1}^{(1)}+ \alpha_{2,2}^{(2)} b_{n-1}^{(1)} + &\cdots &+& \alpha_{2,1}^{(s)} a_{n-1}^{(s-1)}+ \alpha_{2,2}^{(s)} b_{n-1}^{(s-1)}, \\ 
				a_{n}^{(1)} & = &  a_{n-1}, \\
				b_{n}^{(1)} & = &  b_{n-1},\\
				&\vdots&\\
				a_{n}^{(s-1)} & = &  a_{n-s+1}, \\
				b_{n}^{(s-1)} & = &  b_{n-s+1}.      
			\end{array}
		\end{equation*}
		Using \cite[Lemma 2.1]{NSz_Power} (or \cite{NL_grow})
		the matrix form of the system above is 
		\begin{equation}\label{eqsys:ext_vec}
			\mathbf{w}_n={\mathbf{M}^{(s)}}\mathbf{w}_{n-1}={(\mathbf{M}^{(s)})}^n\mathbf{w}_{0},
		\end{equation}
		where for $i\geq s-1$ 
		\[\mathbf{w}_i=
		\begin{pmatrix} a_{i}^{(0)}, b_{i}^{(0)}, a_{i}^{(1)}, b_{i}^{(1)}, \ldots, a_{i}^{(s-1)},b_{i}^{(s-1)}  \end{pmatrix}^T = 
		\begin{pmatrix} a_{i}, b_{i}, a_{i-1}, b_{i-1}, \ldots, a_{i-s+1}, b_{i-s+1}  \end{pmatrix}^T\!.\]
		The initial vector is obviously \[\mathbf{w}_0=\begin{pmatrix} a_{s-1}, b_{s-1}, a_{s-2}, b_{s-2}, \ldots, a_{0}, b_{0}  \end{pmatrix}^T\!,\]
		furthermore
		\[\mathbf{M}^{(s)}=\begin{pmatrix}  
			\alpha_{1,1}^{(1)} & \alpha_{1,2}^{(1)} & \alpha_{1,1}^{(2)} & \alpha_{1,2}^{(2)} &\cdots & \alpha_{1,1}^{(s-1)} & \alpha_{1,2}^{(s-1)} & \alpha_{1,1}^{(s)} & \alpha_{1,2}^{(s)}\\
			\alpha_{2,1}^{(1)} & \alpha_{2,2}^{(1)} & \alpha_{2,1}^{(2)} & \alpha_{2,2}^{(2)} &\cdots & \alpha_{2,1}^{(s-1)} & \alpha_{2,2}^{(s-1)}  & \alpha_{2,1}^{(s)} & \alpha_{2,2}^{(s)}\\
			1 &   &   &   &  &   &   &   &   \\
			& 1 &   &   &   &   &   &   &   \\
			&   & 1 &   &  &   &   &   &   \\
			&   &   & 1 &   &   &   &   &   \\
			&   &   &   & \ddots\\
			&   &   &   &  & 1  &   &   &   \\
			&   &   &   &  &   &  1 &   &  
		\end{pmatrix}_{\!\!2s\times 2s}=
		\left|\arraycolsep=1pt \begin{array}{cccc:c} 
			\mathbf{A}_1  & \mathbf{A}_2  &\cdots & \mathbf{A}_{s-1} & \mathbf{A}_s\\\hdashline
			\\
			&   & \mathbf{I}   &   & \mathbf{0}    \\
		\end{array} \right|,\]
		a block matrix with suitable unit matrix $\mathbf{I}$ and zero matrix $\mathbf{0}$.  (We do not indicate usually the zero entries of the matrices.)  
		
		Applying Lemma 6 of \cite{NL_grow} again, the characteristic polynomial of the coefficient matrix $\mathbf{M}^{(s)}$ and the characteristic polynomial of any component sequence coincide.
		Consequently, we have to determine the characteristic polynomial 
		\[p^{(s)}(x)=|\mathbf{M}^{(s)}-x\mathbf{I}|\] 
		of  $\mathbf{M}^{(s)}$, which 
		yields the common recurrence relation for the sequences $(a_n)$ and $(b_n)$.
		Thus
		\begin{equation}\label{eq:pk_def}
			p^{(s)}(x) =\left|\begin{array}{ccccccc:cc}
				\alpha_{1,1}^{(1)}-x & \alpha_{1,2}^{(1)} & \alpha_{1,1}^{(2)} & \alpha_{1,2}^{(2)} &\cdots & \alpha_{1,1}^{(s-1)} & \alpha_{1,2}^{(s-1)} & \alpha_{1,1}^{(s)} & \alpha_{1,2}^{(s)}\\
				\alpha_{2,1}^{(1)} & \alpha_{2,2}^{(1)}-x & \alpha_{2,1}^{(2)} & \alpha_{2,2}^{(2)} &\cdots & \alpha_{2,1}^{(s-1)} & \alpha_{2,2}^{(s-1)}  & \alpha_{2,1}^{(s)} & \alpha_{2,2}^{(s)}\\
				1 &   & -x  &   &  &   &   &   &   \\
				& 1 &   & -x  &   &   &   &   &   \\
				&   & 1 &   &  &   &   &   &   \\
				&   &   & 1 &   &   &   &   &   \\
				&   &   &   & \ddots &   & \ddots\\ \hdashline
				&   &   &   &  & 1  &   & -x  &   \\
				&   &   &   &  &   &  1 &   & -x 
			\end{array} \right|.
		\end{equation}
		If $s=1$, then we can easily have
		\begin{equation*}\label{eq:pk=1}
			\begin{split}
				p^{(1)}(x) &=\left|\begin{array}{cc}
					\alpha_{1,1}^{(1)}-x & \alpha_{1,2}^{(1)} \\
					\alpha_{2,1}^{(1)} & \alpha_{2,2}^{(1)}-x 
				\end{array} \right|
				= x^2- \left(\alpha_{1,1}^{(1)}+\alpha_{2,2}^{(1)}\right)x +
				\alpha_{1,1}^{(1)}\alpha_{2,2}^{(1)} - \alpha_{1,2}^{(1)}\alpha_{2,1}^{(1)} \\
				&= x^2-\tr(\mathbf{A}_1)x +|\mathbf{A}_1|,
			\end{split}
		\end{equation*}
		hence \[\mathbf{c}^{(1)}=\begin{pmatrix} \tr(\mathbf{A}_1)\\ -|\mathbf{A}_1| \end{pmatrix}.\]    
		
		Now we prove \eqref{eqvect_recu}  by induction on $s\geq2$.
		If $s=2$, then  
		\begin{equation*}\label{eq:pk=2}
			\begin{split}
				p^{(2)}(x) &=\left|\begin{array}{ccll}
					\alpha_{1,1}^{(1)}-x & \alpha_{1,2}^{(1)} & \alpha_{1,1}^{(2)} & \alpha_{1,2}^{(2)} \\
					\alpha_{2,1}^{(1)} & \alpha_{2,2}^{(1)}-x & \alpha_{2,1}^{(1)} & \alpha_{2,2}^{(1)}\\
					1 &0 &-x & 0\\
					0&1 &0 & -x\\
				\end{array} \right| \\
				&= x^4-\tr(\mathbf{A}_1)x^3 +(|\mathbf{A}_1|- \tr(\mathbf{A}_2))x^2 +\Det{\mathbf{A}^{(1,2)}}x+ |\mathbf{A}_2|
			\end{split}
		\end{equation*}
		with
		\[\arraycolsep=-3pt \mathbf{c}^{(2)}=
		\begin{pmatrix} 
			\tr(\mathbf{A}_1)&\\ 
			-|\mathbf{A}_1| &  \ \  + \tr(\mathbf{A}_2)\\ 
			& -\Det{\mathbf{A}^{(1,2)}}\\
			& -|\mathbf{A}_2|  
		\end{pmatrix}=  
		\begin{pmatrix}  \mathbf{c}^{(1)}\\   \\  0\\   0    \end{pmatrix}+
		\begin{pmatrix} 
			0\\ 
			\tr(\mathbf{A}_2)\\ 
			-\Det{\mathbf{A}^{(1,2)}}\\
			0 
		\end{pmatrix}+
		\begin{pmatrix} 
			0\\ 0\\ 0\\
			-|\mathbf{A}_2|  
		\end{pmatrix}.\]
		Moreover, a straightforward calculation shows that   
		\[\arraycolsep=-3pt \mathbf{c}^{(3)}=
		\begin{pmatrix} 
			\tr(\mathbf{A}_1)& &\\ 
			-|\mathbf{A}_1|&  \ \  + \tr(\mathbf{A}_2) &\\
			& -\Det{\mathbf{A}^{(1,2)}} &  \ \  + \tr(\mathbf{A}_3)\\
			& -|\mathbf{A}_2|& -\Det{\mathbf{A}^{(1,3)}} \\
			& & -\Det{\mathbf{A}^{(2,3)}} \\
			& & -|\mathbf{A}_3|
		\end{pmatrix}=
		\begin{pmatrix}\\ \mathbf{c}^{(2)}\\ \\ \\0\\0\end{pmatrix}+
		\begin{pmatrix}0\\ 0\\
			\tr(\mathbf{A}_3)\\
			-\Det{\mathbf{A}^{(1,3)}}\\
			-\Det{\mathbf{A}^{(2,3)}} \\
			0
		\end{pmatrix}+
		\begin{pmatrix}0\\ 0\\  0\\  0\\  0 \\  -|\mathbf{A}_3|
		\end{pmatrix}\]
	holds for $k=3$.
		
		Suppose now that \eqref{eqvect_recu} is true for up to $s-1$. 
		Expand the determinant of \eqref{eq:pk_def} along the last two rows consecutively. 
		We distinguish three cases.
		\begin{description}
			\item[Case 1] 
			If we choose $-x$ from both of the last two rows, then we obtain $(-x)^2|\mathbf{M}^{(s-1)}|=x^2\cdot p^{(s-1)}(x)$. This results in the first vector in the sum \eqref{eqvect_recu}.
			
			\item[Case 2]   If we choose twice 1, then we have to choose 1 again from the rows $2s-2,2s-3,\dots,4,3$ step by step (to avoid the product becoming zero), and we have $|\mathbf{A}_k|$. This yields the last vector of \eqref{eqvect_recu}.
			
			\item[Case 3]   Now, we choose exactly one $-x$ entry from the last two rows. 
		
			Firstly, we chose $-x$ from the last row and $1$ from the penultimate row, as Matrix~1 shows in Figure~\ref{fig:Matrix_1}. We denoted the chosen items with rounded squares.
			To proceed, we have two possibilities to find a non-zero entry in row $(2s-2)$. If we have $-x$, then the next non-zero element is necessarily $1$ from row $(2s-3)$. Similarly, in row $(2s-4)$, we have again two choices, and so on. 
			
			Assume that the pair $(-x,1)$ was chosen consecutively $(s-\ell)$ times in total, where $0\leq \ell \leq s-1$. Then, apart from rows 2 and 1, we continue by opting for only value 1. Note that a double 1 choice predestinates the ending in order to obtain a non-zero product (see Case 2).
			This is the green minor matrix in Figure~\ref{fig:Matrix_1} having determinant 1. Lastly, from the first two rows, the minor matrix is 
			\[\text{ either }
			\left(\begin{array}{cc}
				\alpha_{1,2}^{(\ell)}&	\alpha_{1,1}^{(s)}  \\
				\alpha_{2,2}^{(\ell)}&	\alpha_{2,1}^{(s)} 
			\end{array}\right) 
			\text{ or }
			\left(\begin{array}{cc}
				\alpha_{1,2}^{(1)}&	\alpha_{1,1}^{(s)}  \\
				\alpha_{2,2}^{(1)}-x&	\alpha_{2,1}^{(s)}
			\end{array}\right) .
			\]
			(The second matrix occurs when no double 1 appears in the expansion only pairs $(-x,1)$.)
				\begin{figure} [ht!]
				\centering
				\[
				\setlength{\arraycolsep}{1.5pt}
				\NiceMatrixOptions{code-for-first-col=\color{blue} \scriptstyle, cell-space-limits=0pt}
				\begin{pNiceArray}{*{15}{w{c}{0.6cm}}:*{2}{w{c}{0.6cm}}}[first-col] 
					1&	\ast\!-\!x& \ast & \ast & \ast &\cdots & \ast& \ast&\alpha_{1,1}^{(\ell)} & \Block[draw,fill=blue!15,rounded-corners]{2-1}{}\alpha_{1,2}^{(\ell)} &\ast&\cdots &\ast&\ast & \ast & \ast & \Block[draw,fill=blue!15,rounded-corners]{2-1}{} \alpha_{1,1}^{(s)} & \alpha_{1,2}^{(s)}\\
					2&	\ast & \ast\!-\!x & \ast & \ast &\cdots & \ast&\ast& \alpha_{2,1}^{(\ell)} & \alpha_{2,2}^{(\ell)} &\ast&\cdots&\ast&\ast & \ast & \ast  & \alpha_{2,1}^{(s)} & \alpha_{2,2}^{(s)}\\
					3	   &\Block[draw,fill=green!15,rounded-corners]{7-7}{}  	\Block[draw,fill=blue!15,rounded-corners]{1-1}{}1 && -x\\
					\Vdots && \Block[draw,fill=blue!15,rounded-corners]{1-1}{}1 && -x\\
					&.&& \Block[draw,fill=blue!15,rounded-corners]{1-1}{}1 && \ddots\\
					&&\Ddots^{\color{blue} 2(\ell-1)  \text{ times choose } 1}&& \Block[draw,fill=blue!15,rounded-corners]{1-1}{}1 && -x\\
					&&&&& \ddots && -x\\
					2\ell-1&&&&&& \Block[draw,fill=blue!15,rounded-corners]{1-1}{}1 && -x\\
					2\ell&&&&&.&& \Block[draw,fill=blue!15,rounded-corners]{1-1}{}1 && -x\\
					2\ell+1&&&&&&&& \Block[draw,fill=blue!15,rounded-corners]{1-1}{}1 && -x\\
					2\ell+2&&&&&&&.&& 1 && \ddots &. \\ 
					\Vdots &&&&&&&&\Ddots^{\color{blue}(s-\ell) \text{ times choose } 1}&& \Block[draw,fill=blue!15,rounded-corners]{1-1}{}1  && -x&\Ddots^{\color{blue}(s-\ell) \text{ times choose } -x}\\ 
					2s-4   &&&&&&&&&&& \ddots && \Block[draw,fill=blue!15,rounded-corners]{1-1}{}-x\\ 
					2s-3   &&&&&&&&&&&& \Block[draw,fill=blue!15,rounded-corners]{1-1}{}1  && -x\\ 
					2s-2   &&&&&&&&&&&&& 1 &&\Block[draw,fill=blue!15,rounded-corners]{1-1}{} -x\\ \hdashline
					2s-1   &&&&&&&&&&&&&& \Block[draw,fill=red!15,rounded-corners]{1-1}{} 1  && -x&.\\
					2s     &&&&&&&&&&&&&.&& 1 && \Block[draw,fill=red!15,rounded-corners]{1-1}{}-x 
				\end{pNiceArray}
				\]
				\caption[]{Matrix 1 for calculate the determinant}
				\label{fig:Matrix_1}
			\end{figure}
		
			Summarizing this particular way (the first subcase of Case 3) we have 
			\begin{align*}
				&\sum_{\ell=2}^{s-1} (-x)^{s-\ell} (-1)^{(s-\ell)-1}
				\left|\begin{array}{cc}
					\alpha_{1,2}^{(\ell)}&	\alpha_{1,1}^{(s)}  \\
					\alpha_{2,2}^{(\ell)}&	\alpha_{2,1}^{(s)} 
				\end{array}\right|
				+ (-x)^{s-1} (-1)^{(s-1)-1}
				\left|\begin{array}{cc}
					\alpha_{1,2}^{(1)}&	\alpha_{1,1}^{(s)}  \\
					\alpha_{2,2}^{(1)}-x&	\alpha_{2,1}^{(s)} 
				\end{array}\right| =\\
			&\sum_{\ell=2}^{s-1} x^{s-\ell}
				\left|\begin{array}{cc}
					\alpha_{1,1}^{(s)}&	\alpha_{1,2}^{(\ell)}  \\
					\alpha_{2,1}^{(s)}&	\alpha_{2,2}^{(\ell)} 
				\end{array}\right|
				-x^{s-1} 
				\left(-
				\left|\begin{array}{cc}
					\alpha_{1,1}^{(s)}&	\alpha_{1,2}^{(1)}  \\
					\alpha_{2,1}^{(s)}&	\alpha_{2,2}^{(1)} 
				\end{array}\right|
				+x \alpha_{1,1}^{(s)} 
				\right) =\\
			& \sum_{\ell=1}^{s-1} x^{s-\ell}
			\left| \mathbf{A}_{s,\ell} \right|
			-x^{s}  \alpha_{1,1}^{(s)}.
			\end{align*}
						
			Secondly, we switch the roles of $-x$ and 1, i.e., first choose $1$ from the last row and then $-x$ from row $(2s-1)$ as Matrix~2 shows in Figure~\ref{fig:Matrix_2}. Analogously to the previous description, we keep choosing the pair $(1,-x)$ consecutively for a while, and then proceed with the remaining  entries 1. Then finish the procedure with the first two rows. Finally, we find 
			
		\begin{align*}
			&\sum_{\ell=2}^{s-1}(-x)^{s-\ell} (-1)^{(s-\ell)}
			\left|\begin{array}{cc}
				\alpha_{1,1}^{(\ell)}&	\alpha_{1,2}^{(s)}  \\
				\alpha_{2,1}^{(\ell)}&	\alpha_{2,2}^{(s)} 
			\end{array}\right|
			+(-x)^{s-\ell} (-1)^{(s-\ell)}
			\left|\begin{array}{cc}
				\alpha_{1,1}^{(1)}-x&	\alpha_{1,2}^{(s)}  \\
				\alpha_{2,1}^{(1)}&	\alpha_{2,2}^{(s)} 
			\end{array}\right| =\\
		&\sum_{\ell=2}^{s-1} x^{s-\ell}
		\left|\begin{array}{cc}
			\alpha_{1,1}^{(\ell)}&	\alpha_{1,2}^{(s)}  \\
			\alpha_{2,1}^{(\ell)}&	\alpha_{2,2}^{(s)}
			\end{array}\right|
		+x^{s-1} 
		\left(
		\left|\begin{array}{cc}
			\alpha_{1,1}^{(1)}&	\alpha_{1,2}^{(s)}  \\
			\alpha_{2,1}^{(1)}&	\alpha_{2,2}^{(s)} 
		\end{array}\right|
		-x \alpha_{1,1}^{(s)} 
		\right) =\\
		&\sum_{\ell=1}^{s-1} x^{s-\ell}
		\left| \mathbf{A}_{\ell,s} \right|
		-x^{s}  \alpha_{2,2}^{(s)}. 
	\end{align*}
	
		\begin{figure} [h!]
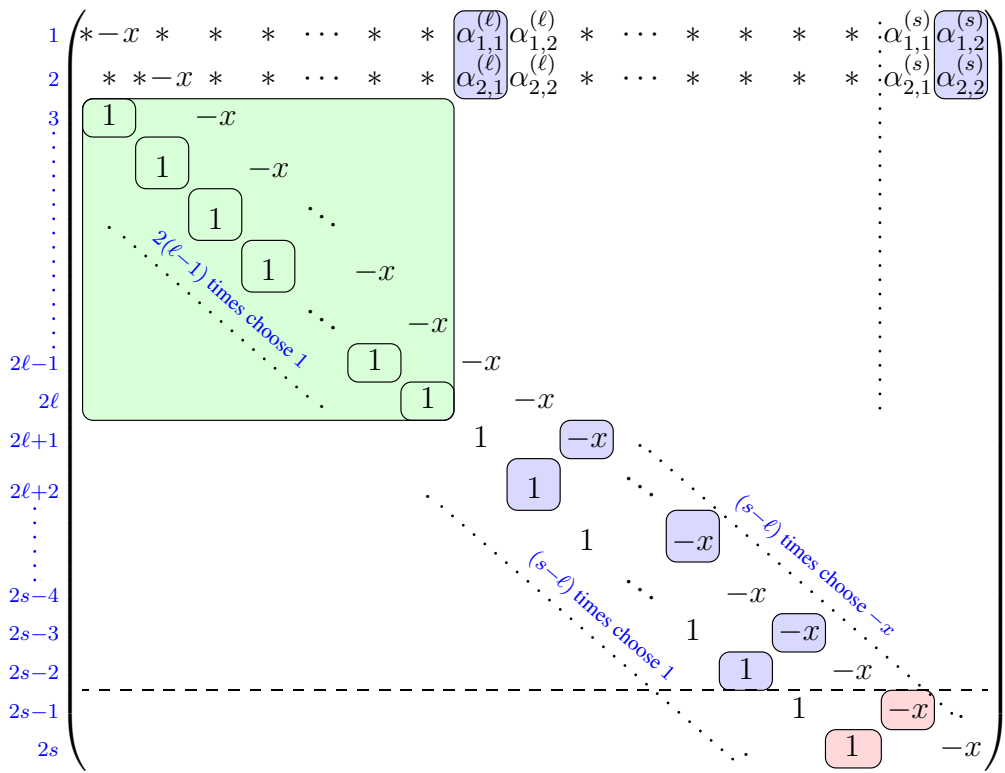

		\centering
		\[
		\setlength{\arraycolsep}{1.5pt}
		\NiceMatrixOptions{code-for-first-col=\color{blue} \scriptstyle, cell-space-limits=0pt}
		\begin{pNiceArray}{*{15}{w{c}{0.6cm}}:*{2}{w{c}{0.6cm}}}[first-col] 
			1&	\ast\!-\!x& \ast & \ast & \ast &\cdots & \ast& \ast& \Block[draw,fill=blue!15,rounded-corners]{2-1}{} 
			\alpha_{1,1}^{(\ell)} & \alpha_{1,2}^{(\ell)} &\ast&\cdots &\ast&\ast & \ast & \ast &  \alpha_{1,1}^{(s)} & \Block[draw,fill=blue!15,rounded-corners]{2-1}{}\alpha_{1,2}^{(s)}\\
			2&	\ast & \ast\!-\!x & \ast & \ast &\cdots & \ast&\ast& \alpha_{2,1}^{(\ell)} & \alpha_{2,2}^{(\ell)} &\ast&\cdots&\ast&\ast & \ast & \ast  & \alpha_{2,1}^{(s)} & \alpha_{2,2}^{(s)}\\
			3	   &\Block[draw,fill=green!15,rounded-corners]{7-7}{}  	\Block[draw,fill=blue!15,rounded-corners]{1-1}{}1 && -x\\
			\Vdots && \Block[draw,fill=blue!15,rounded-corners]{1-1}{}1 && -x\\
			&.&& \Block[draw,fill=blue!15,rounded-corners]{1-1}{}1 && \ddots\\
			&&\Ddots^{\color{blue}2(\ell-1) \text{ times choose } 1}&& \Block[draw,fill=blue!15,rounded-corners]{1-1}{}1 && -x\\
			&&&&& \ddots && -x\\
			2\ell-1&&&&&& \Block[draw,fill=blue!15,rounded-corners]{1-1}{}1 && -x\\
			2\ell&&&&&.&& \Block[draw,fill=blue!15,rounded-corners]{1-1}{}1 && -x\\
			2\ell+1&&&&&&&& 1 && \Block[draw,fill=blue!15,rounded-corners]{1-1}{}-x&.\\
			2\ell+2  &&&&&&&.&& \Block[draw,fill=blue!15,rounded-corners]{1-1}{}1 && \ddots  \\ 
			\Vdots &&&&&&&&\Ddots^{\color{blue}(s-\ell) \text{ times choose } 1}&& 1  &&\Block[draw,fill=blue!15,rounded-corners]{1-1}{} -x&\Ddots^{\color{blue}(s-\ell) \text{ times choose } -x}\\ 
			2s-4   &&&&&&&&&&& \ddots && -x\\ 
			2s-3   &&&&&&&&&&&& 1  && \Block[draw,fill=blue!15,rounded-corners]{1-1}{} -x\\ 
			2s-2   &&&&&&&&&&&&& \Block[draw,fill=blue!15,rounded-corners]{1-1}{}1 &&-x\\ \hdashline
			2s-1   &&&&&&&&&&&&&&  1  && \Block[draw,fill=red!15,rounded-corners]{1-1}{}-x&.\\
			2s     &&&&&&&&&&&&&.&& \Block[draw,fill=red!15,rounded-corners]{1-1}{}1 && -x 
		\end{pNiceArray}
		\]
		\caption[]{Matrix 2 for calculate the determinant}
		\label{fig:Matrix_2}
	\end{figure}
			
		Hence, Case 3 finally returns  with 
		\[	
		\sum_{\ell=1}^{s-1} (x)^{s-\ell} \Det{\mathbf{A}^{(\ell,s)}} - x^{s} \tr{\mathbf{A}_s},
		\]	
		 which gives the second vector in \eqref{eqvect_recu}. 
	\end{description}
	Now the proof is now complete.
\end{proof}

	\begin{remark}
		Theorem \ref{th:main} makes it possible to give the vector  $\mathbf{c}^{(s)}$ explicitly as follows: 
		\[\mathbf{c}^{(k)}=
		\begin{pmatrix}  
			\tr(\mathbf{A}_1)\\ 
			-|\mathbf{A}_1| + \tr(\mathbf{A}_2) \\
			-\Det{\mathbf{A}^{(1,2)}} + \tr(\mathbf{A}_3)\\
			-|\mathbf{A}_2| -\Det{\mathbf{A}^{(1,3)}} + \tr(\mathbf{A}_4)\\
			-\Det{\mathbf{A}^{(2,3)}}-\Det{\mathbf{A}^{(1,4)}}+ \tr(\mathbf{A}_5)\\
			-|\mathbf{A}_3|-\Det{\mathbf{A}^{(2,4)}}-\Det{\mathbf{A}^{(1,5)}}+ \tr(\mathbf{A}_6)\\
			-\Det{\mathbf{A}^{(3,4)}}-\Det{\mathbf{A}^{(2,5)}}-\Det{\mathbf{A}^{(1,6)}}+ \tr(\mathbf{A}_7)\\
			-|\mathbf{A}_4|-\Det{\mathbf{A}^{(3,5)}}-\Det{\mathbf{A}^{(2,6)}}-\Det{\mathbf{A}^{(1,7)}}+ \tr(\mathbf{A}_8)\\
			-\Det{\mathbf{A}^{(4,5)}}-\Det{\mathbf{A}^{(3,6)}}-\Det{\mathbf{A}^{(2,7)}}-\Det{\mathbf{A}^{(1,8)}} + \tr(\mathbf{A}_9)\\ 
			-|\mathbf{A}_5|-\Det{\mathbf{A}^{(4,6)}}-\Det{\mathbf{A}^{(3,7)}}-\cdots -\Det{\mathbf{A}^{(1,9)}} + \tr(\mathbf{A}_{10})\\
			\vdots\\
			-|\mathbf{A}_{\lfloor (s-1)/2\rfloor}- \cdots     +\tr(\mathbf{A}_{s-2})\\
			-\Det{\mathbf{A}^{\lfloor (s-1)/2\rfloor,s-1-\lfloor (s-1)/2\rfloor)}}- \cdots -\Det{\mathbf{A}^{(1,s-2)}}    +\tr(\mathbf{A}_{s-1})\\
			-|\mathbf{A}_{\lfloor s/2\rfloor}|-\Det{\mathbf{A}^{(-1+s/2,1+s/2)}}- \cdots -\Det{\mathbf{A}^{(2,s-2)}} -\Det{\mathbf{A}^{(1,s-1)}}    +\tr(\mathbf{A}_s)\\
			\cdots -\Det{\mathbf{A}^{(3,s-2)}} -\Det{\mathbf{A}^{(2,s-1)}} -\Det{\mathbf{A}^{(1,s)}} \\
			\cdots -\Det{\mathbf{A}^{(4,s-2)}} -\Det{\mathbf{A}^{(3,s-1)}} -\Det{\mathbf{A}^{(2,s)}} \\
			\cdots -\Det{\mathbf{A}^{(5,s-2)}} -\Det{\mathbf{A}^{(4,s-1)}} -\Det{\mathbf{A}^{(3,s)}} \\
			\cdots -\Det{\mathbf{A}^{(6,s-2)}} -\Det{\mathbf{A}^{(5,s-1)}} -\Det{\mathbf{A}^{(4,s)}} \\
			\cdots -\Det{\mathbf{A}^{(7,s-2)}} -\Det{\mathbf{A}^{(6,s-1)}} -\Det{\mathbf{A}^{(5,s)}} \\
			\vdots\\
			-|\mathbf{A}_{s-2}| -\Det{\mathbf{A}^{(s-4,s-2)}} -  \Det{\mathbf{A}^{(s-5,s-1)}} -\Det{\mathbf{A}^{(s-6,s)}} \\
			-\Det{\mathbf{A}^{(s-3,s-2)}} -\Det{\mathbf{A}^{(s-4,s-1)}} -\Det{\mathbf{A}^{(s-5,s)}} \\
			-|\mathbf{A}_{s-2}| -\Det{\mathbf{A}^{(s-3,s-1)}} -\Det{\mathbf{A}^{(s-4,s)}} \\
			-\Det{\mathbf{A}^{(s-2,s-1)}} -\Det{\mathbf{A}^{(s-3,s)}} \\
			-|\mathbf{A}_{s-1}| -\Det{\mathbf{A}^{(s-2,s)}} \\
			-\Det{\mathbf{A}^{(s-1,s)}} \\
			-|\mathbf{A}_s|
		\end{pmatrix}.
		\]
	\end{remark}

	\section{Example}

	The binary recurrence $(u_n)$ is defined by
	\begin{equation} \label{eq:u}
		u_{n}=a u_{n-1}+b u_{n-2},\quad (n\geq 2)
	\end{equation}
	with initial values $u_0=1$, $u_1=a$ (or $u_{-1}=0$).
	The term $u_n$ is interpreted in the work \cite{Ben} as the number of ways to tile a $(1 \times n)$-board (or square grid) using squares with $a$ different colors  and dominoes with $b$ colors. Obviously, if $a=b=1$ then $u_n=F_{n+1}$, the shifted Fibonacci sequence.

	\subsection{Tiling with black and white dominoes}
	
	Now we present a variant of the problem above.
	Fix a positive integer $k\ge1$ as the maximal size of the dominoes we will use. Then consider a $(1\times n)$-board  and tile it from left to right with black and white dominoes such that black domino never follows white domino. Let $t_{n,k}$ denote the number of such tiles. 
	
\begin{figure}[ht!]
		\centering
		\includegraphics[scale=1]{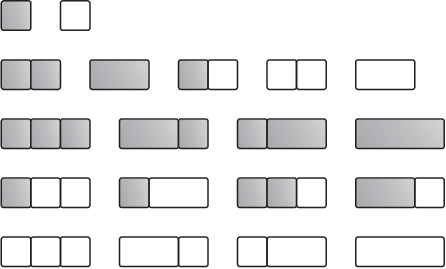} 
		\caption{Tiling with black and white dominoes if $n\leq3$ and $k\geq3$.}
		\label{fig:recu_sys_n3_k3}
	\end{figure}

	Under the conditions above, let $a_{n,k}$ and $b_{n,k}$ denote the number of different tilings ending with black and white dominos, respectively. Now $t_{n,k}=a_{n,k}+b_{n,k}$. The tiles in Figure~\ref{fig:recu_sys_n3_k3} shows the possible cases when $n=1, 2, 3$ and $k\geq3$. It admits $a_{1,k}=1$, $b_{1,k}=1$, $a_{2,k}=2$, $b_{2,k}=3$, $a_{3,k}=4$, $b_{3,k}=8$. Assume that  we let $a_{0,k}=1$ and $b_{1,k}=0$ for any $k$.
	
	It is easy to see by the constraint that 
	\begin{equation}\label{eqsys:example}
		\begin{array}{lclclclcl} 
			a_{n,k} & = &  a_{n-1,k} &+&  a_{n-2,k} &+& \cdots &+&  a_{n-k,k},\\ 
			b_{n,k} & = & a_{n-1,k}+  b_{n-1,k} &+&  a_{n-2,k}+  b_{n-2,k} &+ &\cdots &+&  a_{n-k,k}+  b_{n-k,k},
		\end{array} 
	\end{equation} 
	a specific example of the system \eqref{eqsys:def}.
	Here, for $i,j\geq1$ we have $\mathbf{A}_j=\begin{pmatrix}  1 &0\\  1 & 1 \end{pmatrix}$, further $|\mathbf{A}_j|=1$,  $\tr(\mathbf{A}_j)=2$, and $\Det{\mathbf{A}^{(i,j)}}=2$. Moreover, we find for the coefficients of the recurrence sequences that   
	\[\mathbf{c}^{(1)}=\begin{pmatrix} \tr(\mathbf{A}_1)\\ -|\mathbf{A}_1| \end{pmatrix}= \begin{pmatrix} 2\\ -1 \end{pmatrix},  \quad
	\arraycolsep=-3pt \mathbf{c}^{(2)}=
	\begin{pmatrix} 
		\tr(\mathbf{A}_1)&\\ 
		-|\mathbf{A}_1|&\, + \tr(\mathbf{A}_2)\\ 
		& -\Det{\mathbf{A}^{(1,2)}}\\
		& -|\mathbf{A}_2|  
	\end{pmatrix} =
	\begin{pmatrix} 
		2\\ 
		1\\ 
		-2\\
		-1  
	\end{pmatrix}, \]
	and generally 
	\begin{equation}\label{eqvect_recu_example}
		\mathbf{c}^{(k)}=
		\begin{pmatrix}\\ \\ \\ \\ \mathbf{c}^{(k-1)}\\ \\ \\0\\0\end{pmatrix}+
		\begin{pmatrix} 0\\ \vdots\\ 0\\ 2\\ -2 \\ 	\vdots\\-2 \\ -1
		\end{pmatrix}.
	\end{equation}  
	Writing the vectors $(\mathbf{c}^{(k)})^T$ in triangle form we have the coefficient triangle $\mathcal{C}$ (see Table~\ref{fig:Triangle_coef}), where row $k$ contains the components of $\mathbf{c}^{(k)}$. Thus, the common recurrence relation of sequences $(a_{n,k})$, $(b_{n,k})$, and $(t_{n,k})$  for $k=1,2,3$ (in form \eqref{eq:z-K_vec}) are 
	\[
	\begin{aligned}
		\text{if } k=1, \quad& z_{n,1} =2 z_{n-1,1}-z_{n-2,1} & (n\geq2), \\
		\text{if } k=2, \quad& z_{n,2} =2 z_{n-1,2}+z_{n-2,2} -2z_{n-3,2} -z_{n-4,2} , & (n\geq4), \\
		\text{if } k=3, \quad& z_{n,2} =2 z_{n-1,2}+z_{n-2,2} -3z_{n-4,2} -2z_{n-5,2}-z_{n-6,2}  & (n\geq6).
	\end{aligned}
	\]
	The initial values are 
	$a_{0,k}=1$, $b_{0,k}=0$, $t_{0,k}=1$, 
	$a_{1,k}=1$, $b_{1,k}=1$, $t_{1,k}=2$, and  for $2\leq n\leq k-1$,
	$a_{n,k}=a_{n,k-1}$, $b_{n,k}=b_{n,k-1}$, $t_{n,k}=t_{n,k-1}$.

	\begin{table}[!ht]
		\centering
		\scalebox{0.99}{ \begin{tikzpicture}[->,xscale=0.75,yscale=0.3, auto,swap]
				
				\node (a01) at (-0.5,-1)   {2};
				\node (a02) at (0.5,-1)   {$-1$};

				\node (a02) at (-1.5,-3)   {2};		
				\node (a01) at (-0.5,-3)   {1};
				\node (a01) at (0.5,-3)   {$-2$};
				\node (a02) at (1.5,-3)   {$-1$};

				\node (a01) at (-2.5,-5)   {2};		
				\node (a01) at (-1.5,-5)   {1};
				\node (a02) at (-0.5,-5)   {0};
				\node (a02) at (0.5,-5)   {$-3$};
				\node (a01) at (1.5,-5)   {$-2$};
				\node (a01) at (2.5,-5)   {$-1$};

				\node (a01) at (-3.5,-7)   {2};
				\node (a01) at (-2.5,-7)   {1};		
				\node (a01) at (-1.5,-7)   {0};
				\node (a02) at (-0.5,-7)   {$-1$};
				\node (a02) at (0.5,-7)   {$-4$};
				\node (a01) at (1.5,-7)   {$-3$};
				\node (a01) at (2.5,-7)   {$-2$};
				\node (a01) at (3.5,-7)   {$-1$};
				
				\node (a01) at (-4.5,-9)   {2};
				\node (a01) at (-3.5,-9)   {1};			
				\node (a01) at (-2.5,-9)   {0};		
				\node (a01) at (-1.5,-9)   {$-1$};
				\node (a02) at (-0.5,-9)   {$-2$};
				\node (a02) at (0.5,-9)   {$-5$};
				\node (a01) at (1.5,-9)   {$-4$};
				\node (a01) at (2.5,-9)   {$-3$};
				\node (a01) at (3.5,-9)   {$-2$};
				\node (a01) at (4.5,-9)   {$-1$};
				
				\node (a01) at (-5.5,-11)   {2};
				\node (a01) at (-4.5,-11)   {1};	
				\node (a01) at (-3.5,-11)   {0};			
				\node (a01) at (-2.5,-11)   {$-1$};		
				\node (a01) at (-1.5,-11)   {$-2$};
				\node (a02) at (-0.5,-11)   {$-3$};
				\node (a02) at (0.5,-11)   {$-6$};
				\node (a01) at (1.5,-11)   {$-5$};
				\node (a01) at (2.5,-11)   {$-4$};
				\node (a01) at (3.5,-11)   {$-3$};
				\node (a01) at (4.5,-11)   {$-2$};
				\node (a01) at (5.5,-11)   {$-1$};
				
				\node (a01) at (-6.5,-13)   {2};
				\node (a01) at (-5.5,-13)   {1};	
				\node (a01) at (-4.5,-13)   {0};			
				\node (a01) at (-3.5,-13)   {$-1$};		
				\node (a01) at (-2.5,-13)   {$-2$};
				\node (a02) at (-1.5,-13)   {$-3$};
				\node (a02) at (-0.5,-13)   {$-4$};
				\node (a01) at (0.5,-13)   {$-7$};  
				\node (a02) at (1.5,-13)   {$-6$};
				\node (a01) at (2.5,-13)   {$-5$};
				\node (a01) at (3.5,-13)   {$-4$};
				\node (a01) at (4.5,-13)   {$-3$};
				\node (a01) at (5.5,-13)   {$-2$};
				\node (a01) at (6.5,-13)   {$-1$}; 			
		\end{tikzpicture}}
		\caption{Coefficient triangle $\mathcal{C}$} 
		\label{fig:Triangle_coef}
	\end{table}

	\begin{remark}
		From the first row of the system \eqref{eqsys:example} it can be seen that the sequences $(a_{n,k})$ are the (nonzero) $k$-generalized Fibonacci sequences (or $k$-bonacci sequences), and in case of the first few $k$, they are appear in {The On-Line Encyclopedia of Integer Sequences (OEIS)}~\cite{oeis}, see the sequences A000007, A000012, A000045, A000073, A000078, A001591, A001592, A066178, A122189, A079262.
		Moreover, $(a_{n,0})$, $(a_{n,1})$, $(a_{n,2})$, and $(t_{n,0})$, $(t_{n,1})$, $(t_{n,2})$, $(t_{n,3})$, $(t_{n,4})$ are also in OEIS~\cite{oeis} as sequences A000004, A001477, A023610, and A000007, A000027, A001629, A073778, A118898.
	\end{remark}

	\section{Acknowledgment}
	
	The research was supported in part by National Research, Development and Innovation Office Grant 2019-2.1.11-T\'ET-2020-00165. 
	L.~Szalay was supported by the Hungarian National Foundation for Scientific Research Grant No.~128088, and No.~130909.

	\bigskip
	
	
	\noindent\textbf{Conflict of interest.} The authors declare that they have no conflict of interest.
	\medskip
	
	\noindent\textbf{Data availability.} No data was used for the research described in the article.
	\medskip
	

	\bibliography{Recsys_bib.bib} 
	\bibliographystyle{IEEEtranN}

\end{document}